\documentclass[12pt]{article}
\usepackage{amsmath,amssymb,amsthm,latexsym,amscd}

\title{Combinations of structures\footnote{
\newline\indent \ \ \ The research is partially supported by
Committee of Science in Education and Science Ministry of the
Republic of Kazakhstan, Grant No. 0830/GF4. } }
\author{Sergey V.
Sudoplatov\footnote{sudoplat@math.nsc.ru}}
\date{}
\begin{document}
\maketitle

\begin{abstract}
We investigate combinations of structures by families of
structures relative to families of unary predicates and
equivalence relations. Conditions preserving $\omega$-categoricity
and Ehrenfeuchtness under these combinations are characterized.
The notions of $e$-spectra are introduced and possibilities for
$e$-spectra are described.

\end{abstract}

The aim of the paper is to introduce operators (similar to
\cite{Wo, Su99, Su031, SuCCMCT}) on classes of structures
producing structures approximating given structure, as well as to
study properties of these operators.

In Section 1 we define $P$-operators, $E$-operators, and
corresponding combinations of structures. In Section 2 we
characterize the preservation of $\omega$-categoricity for
$P$-combinations and $E$-combinations as well as Ehrenfeuchtness
for $P$-combinations. In Section 3 we pose and investigate
questions on variations of structures under $P$-operators and
$E$-operators. The notions of $e$-spectra for $P$-operators and
$E$-operators are introduced in Section 4. Here values for
$e$-spectra are described. In Section 5 the preservation of
Ehrenfeuchtness for $E$-combinations is characterized.

Throughout the paper we consider structures of relational
languages.

\section{$P$-operators, $E$-operators, combinations}

Let $P=(P_i)_{i\in I}$, be a family of nonempty unary predicates,
$(\mathcal{A}_i)_{i\in I}$ be a family of structures such that
$P_i$ is the universe of $\mathcal{A}_i$, $i\in I$, and the
symbols $P_i$ are disjoint with languages for the structures
$\mathcal{A}_j$, $j\in I$. The structure
$\mathcal{A}_P\rightleftharpoons\bigcup\limits_{i\in
I}\mathcal{A}_i$\index{$\mathcal{A}_P$} expanded by the predicates
$P_i$ is the {\em $P$-union}\index{$P$-union} of the structures
$\mathcal{A}_i$, and the operator mapping $(\mathcal{A}_i)_{i\in
I}$ to $\mathcal{A}_P$ is the {\em
$P$-operator}\index{$P$-operator}. The structure $\mathcal{A}_P$
is called the {\em $P$-combination}\index{$P$-combination} of the
structures $\mathcal{A}_i$ and denoted by ${\rm
Comb}_P(\mathcal{A}_i)_{i\in I}$\index{${\rm
Comb}_P(\mathcal{A}_i)_{i\in I}$} if
$\mathcal{A}_i=(\mathcal{A}_P\upharpoonright
A_i)\upharpoonright\Sigma(\mathcal{A}_i)$, $i\in I$. Structures
$\mathcal{A}'$, which are elementary equivalent to ${\rm
Comb}_P(\mathcal{A}_i)_{i\in I}$, will be also considered as
$P$-combinations.

By the definition, without loss of generality we can assume for
${\rm Comb}_P(\mathcal{A}_i)_{i\in I}$ that all languages
$\Sigma(\mathcal{A}_i)$ coincide interpreting new predicate
symbols for $\mathcal{A}_i$ by empty relation.

Clearly, all structures $\mathcal{A}'\equiv {\rm
Comb}_P(\mathcal{A}_i)_{i\in I}$ are represented as unions of
their restrictions $\mathcal{A}'_i=(\mathcal{A}'\upharpoonright
P_i)\upharpoonright\Sigma(\mathcal{A}_i)$ if and only if the set
$p_\infty(x)=\{\neg P_i(x)\mid i\in I\}$ is inconsistent. If
$\mathcal{A}'\ne{\rm Comb}_P(\mathcal{A}'_i)_{i\in I}$, we write
$\mathcal{A}'={\rm Comb}_P(\mathcal{A}'_i)_{i\in
I\cup\{\infty\}}$, where
$\mathcal{A}'_\infty=\mathcal{A}'\upharpoonright
\bigcap\limits_{i\in I}\overline{P_i}$, maybe applying
Morleyzation. Moreover, we write ${\rm
Comb}_P(\mathcal{A}_i)_{i\in I\cup\{\infty\}}$\index{${\rm
Comb}_P(\mathcal{A}_i)_{i\in I\cup\{\infty\}}$} for ${\rm
Comb}_P(\mathcal{A}_i)_{i\in I}$ with the empty structure
$\mathcal{A}_\infty$.

Note that if all predicates $P_i$ are disjoint, a structure
$\mathcal{A}_P$ is a $P$-combination and a disjoint union of
structures $\mathcal{A}_i$ \cite{Wo}. In this case the
$P$-combination $\mathcal{A}_P$ is called {\em
disjoint}.\index{$P$-combination!disjoint} Clearly, for any
disjoint $P$-combination $\mathcal{A}_P$, ${\rm
Th}(\mathcal{A}_P)={\rm Th}(\mathcal{A}'_P)$, where
$\mathcal{A}'_P$ is obtained from $\mathcal{A}_P$ replacing
$\mathcal{A}_i$ by pairwise disjoint
$\mathcal{A}'_i\equiv\mathcal{A}_i$, $i\in I$. Thus, in this case,
similar to structures the $P$-operator works for the theories
$T_i={\rm Th}(\mathcal{A}_i)$ producing the theory $T_P={\rm
Th}(\mathcal{A}_P)$\index{$T_P$}, which is denoted by ${\rm
Comb}_P(T_i)_{i\in I}$.\index{${\rm Comb}_P(T_i)_{i\in I}$}

On the opposite side, if all $P_i$ coincide then
$P_i(x)\equiv(x\approx x)$ and removing the symbols $P_i$ we get
the restriction of $\mathcal{A}_P$ which is the combination of the
structures $\mathcal{A}_i$ \cite{Su031, SuCCMCT}.

For an equivalence relation $E$ replacing disjoint predicates
$P_i$ by $E$-classes we get the structure
$\mathcal{A}_E$\index{$\mathcal{A}_E$} being the {\em
$E$-union}\index{$E$-union} of the structures $\mathcal{A}_i$. In
this case the operator mapping $(\mathcal{A}_i)_{i\in I}$ to
$\mathcal{A}_E$ is the {\em $E$-operator}\index{$E$-operator}. The
structure $\mathcal{A}_E$ is also called the {\em
$E$-combination}\index{$E$-combination} of the structures
$\mathcal{A}_i$ and denoted by ${\rm Comb}_E(\mathcal{A}_i)_{i\in
I}$\index{${\rm Comb}_E(\mathcal{A}_i)_{i\in I}$}; here
$\mathcal{A}_i=(\mathcal{A}_E\upharpoonright
A_i)\upharpoonright\Sigma(\mathcal{A}_i)$, $i\in I$. Similar
above, structures $\mathcal{A}'$, which are elementary equivalent
to $\mathcal{A}_E$, are denoted by ${\rm
Comb}_E(\mathcal{A}'_j)_{j\in J}$, where $\mathcal{A}'_j$ are
restrictions of $\mathcal{A}'$ to its $E$-classes.

If $\mathcal{A}_E\prec\mathcal{A}'$, the restriction
$\mathcal{A}'\upharpoonright(A'\setminus A_E)$ is denoted by
$\mathcal{A}'_\infty$. Clearly,
$\mathcal{A}'=\mathcal{A}'_E\coprod\mathcal{A}'_\infty$, where
$\mathcal{A}'_E={\rm Comb}_E(\mathcal{A}'_i)_{i\in I}$,
$\mathcal{A}'_i$ is a restriction of $\mathcal{A}'$ to its
$E$-class containing the universe $A_i$, $i\in I$.

Considering an $E$-combination $\mathcal{A}_E$ we will identify
$E$-classes $A_i$ with structures $\mathcal{A}_i$.

Clearly, the nonempty structure $\mathcal{A}'_\infty$ exists if
and only if $I$ is infinite.

Notice that any $E$-operator can be interpreted as $P$-operator
replacing or naming $E$-classes for $\mathcal{A}_i$ by unary
predicates $P_i$. For infinite $I$, the difference between
`replacing' and `naming' implies that $\mathcal{A}_\infty$ can
have unique or unboundedly many $E$-classes returning to the
$E$-operator.

Thus, for any $E$-combination $\mathcal{A}_E$, ${\rm
Th}(\mathcal{A}_E)={\rm Th}(\mathcal{A}'_{E})$, where
$\mathcal{A}'_{E}$ is obtained from $\mathcal{A}_E$ replacing
$\mathcal{A}_i$ by pairwise disjoint
$\mathcal{A}'_i\equiv\mathcal{A}_i$, $i\in I$. In this case,
similar to structures the $E$-operator works for the theories
$T_i={\rm Th}(\mathcal{A}_i)$ producing the theory $T_E={\rm
Th}(\mathcal{A}_E)$,\index{$T_E$} which is denoted by ${\rm
Comb}_E(T_i)_{i\in I}$\index{${\rm Comb}_E(T_i)_{i\in I}$}, by
$\mathcal{T}_E$\index{$\mathcal{T}_E$}, or by ${\rm
Comb}_E\mathcal{T}$\index{${\rm Comb}_E\mathcal{T}$}, where
$\mathcal{T}=\{T_i\mid i\in I\}$.

Note that $P$-combinations and $E$-unions can be interpreted by
randomizations \cite{AK} of structures.

Sometimes we admit that combinations ${\rm
Comb}_P(\mathcal{A}_i)_{i\in I}$ and ${\rm
Comb}_E(\mathcal{A}_i)_{i\in I}$ are expanded by new relations or
old relations are extended by new tuples. In these cases the
combinations will be denoted by ${\rm
EComb}_P(\mathcal{A}_i)_{i\in I}$ and ${\rm
EComb}_E(\mathcal{A}_i)_{i\in I}$ respectively.

\section{$\omega$-categoricity and Ehrenfeuchtness for combinations}

\medskip
{\bf Proposition 2.1.} {\em If predicates $P_i$ are pairwise
disjoint, the languages $\Sigma(\mathcal{A}_i)$ are at most
countable, $i\in I$, $|I|\leq\omega$, and the structure
$\mathcal{A}_P$ is infinite then the theory ${\rm
Th}(\mathcal{A}_P)$ is $\omega$-categorical if and only if $I$ is
finite and each structure $\mathcal{A}_i$ is either finite or
$\omega$-categorical.}

\medskip
{\bf\em Proof.} If $I$ is infinite or there is an infinite
structure $\mathcal{A}_i$ which is not $\omega$-categorical then
$T={\rm Th}(\mathcal{A}_P)$ has infinitely many $n$-types, where
$n=1$ if $|I|\geq\omega$ and $n=n_0$ for ${\rm Th}(\mathcal{A}_i)$
with infinitely many $n_0$-types. Hence by Ryll-Nardzewski Theorem
${\rm Th}(\mathcal{A}_P)$ is not $\omega$-categorical.

If ${\rm Th}(\mathcal{A}_P)$ is $\omega$-categorical then by
Ryll-Nardzewski Theorem having finitely many $n$-types for each
$n\in\omega$, we have both finitely many predicates $P_i$ and
finitely many $n$-types for each $P_i$-restriction, i.~e., for
${\rm Th}(\mathcal{A}_i)$.~$\Box$

\medskip
Notice that Proposition 2.1 is not true if a $P$-combination is
not disjoint: taking, for instance, a graph $\mathcal{A}_1$ with a
set $P_1$ of vertices and with infinitely many $R_1$-edges such
that all vertices have degree $1$, as well as taking a graph
$\mathcal{A}_2$ with the same set $P_1$ of vertices and with
infinitely many $R_2$-edges such that all vertices have degree
$1$, we can choose edges such that $R_1\cap R_2=\varnothing$, each
vertex in $P_1$ has $(R_1\cup R_2)$-degree $2$, and alternating
$R_1$- and $R_2$-edges there is an infinite sequence of $(R_1\cup
R_2)$-edges. Thus, $\mathcal{A}_1$ and $\mathcal{A}_2$ are
$\omega$-categorical whereas ${\rm
Comb}(\mathcal{A}_1,\mathcal{A}_2)$ is not.

Note also that Proposition 2.1 does not hold replacing
$\mathcal{A}_P$ by $\mathcal{A}_E$. Indeed, taking infinitely many
infinite $E$-classes with structures of the empty languages we get
an $\omega$-categorical structure of the equivalence relation $E$.
At the same time, Proposition 2.1 is preserved if there are
finitely many $E$-classes. In general case $\mathcal{A}_E$ does
not preserve the $\omega$-categoricity if and only if
$E_i$-classes {\em approximate}\index{$E_i$-classes!approximating
infinitely many $n$-types} infinitely many $n$-types for some
$n\in\omega$, i.~e., there are infinitely many $n$-types
$q_m(\bar{x})$, $m\in\omega$, such that for any $m\in\omega$,
$\varphi_j(\bar{x})\in q_j(\bar{x})$, $j\leq m$, and classes
$E_{k_1},\ldots,E_{k_m}$, all formulas $\varphi_j(\bar{x})$ have
realizations in $A_E\setminus\bigcup\limits_{r=1}^m E_{k_r}$.
Indeed, assuming that all $\mathcal{A}_i$ are $\omega$-categorical
we can lose the $\omega$-categoricity for ${\rm
Th}(\mathcal{A}_E)$ only having infinitely many $n$-types (for
some $n$) inside $\mathcal{A}_\infty$. Since all $n$-types  in
$\mathcal{A}_\infty$ are locally (for any formulas in these types)
realized in infinitely many $\mathcal{A}_i$, $E_i$-classes
approximate infinitely many $n$-types and ${\rm
Th}(\mathcal{A}_E)$ is not $\omega$-categorical. Thus, we have the
following

\medskip
{\bf Proposition 2.2.} {\em If the languages
$\Sigma(\mathcal{A}_i)$ are at most countable, $i\in I$,
$|I|\leq\omega$, and the structure $\mathcal{A}_E$ is infinite
then the theory ${\rm Th}(\mathcal{A}_E)$ is $\omega$-categorical
if and only if each structure $\mathcal{A}_i$ is either finite or
$\omega$-categorical, and $I$ is either finite, or infinite and
$E_i$-classes do not approximate infinitely many $n$-types for any
$n\in\omega$.}

\medskip
As usual we denote by $I(T,\lambda)$ the number of pairwise
non-isomorphic models of $T$ having the cardinality $\lambda$.

Recall that a theory $T$ is {\em
Ehrenfeucht}\index{Theory!Ehrenfeucht} if $T$ has finitely many
countable models ($I(T,\omega)<\omega$) but is not
$\omega$-categorical ($I(T,\omega)>1$). A structure with an
Ehrenfeucht theory is also {\em
Ehrenfeucht}\index{Structure!Ehrenfeucht}.

\medskip
{\bf Theorem 2.3.} {\em If predicates $P_i$ are pairwise disjoint,
the languages $\Sigma(\mathcal{A}_i)$ are at most countable, $i\in
I$, and the structure $\mathcal{A}_P$ is infinite then the theory
${\rm Th}(\mathcal{A}_P)$ is Ehrenfeucht if and only if the
following conditions hold:

{\rm (a)} $I$ is finite;

{\rm (b)} each structure $\mathcal{A}_i$ is either finite, or
$\omega$-categorical, or Ehrenfeucht;

{\rm (c)} some $\mathcal{A}_i$ is Ehrenfeucht.}

\medskip
{\bf\em Proof.} If $I$ is finite, each structure $\mathcal{A}_i$
is either finite, or $\omega$-categorical, or Ehrenfeucht, and
some $\mathcal{A}_i$ is Ehrenfeucht then $T={\rm
Th}(\mathcal{A}_P)$ is Ehrenfeucht since each model of $T$ is
composed of disjoint models with universes $P_i$ and
\begin{equation}\label{cs1}
I(T,\omega)=\prod\limits_{i\in I}I({\rm
Th}(\mathcal{A}_i),\min\{|A_i|,\omega\}).
\end{equation}

Now if $I$ is finite and all $\mathcal{A}_i$ are
$\omega$-categorical then by (\ref{cs1}), $I(T,\omega)=1$, and if
some $I({\rm Th}(\mathcal{A}_i),\omega)\geq\omega$ then again by
(\ref{cs1}), $I(T,\omega)\geq\omega$.

Assuming that $|I|\geq\omega$ we have to show that the
non-$\omega$-categorical theory $T$ has infinitely many countable
models. Assuming on contrary that $I(T,\omega)<\omega$, i.~e., $T$
is Ehrenfeucht, we have a nonisolated powerful type $q(\bar{x})\in
S(T)$ \cite{Be}, i.~e., a type such that any model of $T$
realizing $q(\bar{x})$ realizes all types in $S(T)$. By the
construction of disjoint union, $q(\bar{x})$ should have a
realization of the type $p_\infty(x)=\{\neg P_i(x)\mid i\in I\}$.
Moreover, if some ${\rm Th}(\mathcal{A}_i)$ is not
$\omega$-categorical for infinite $A_i$ then $q(\bar{x})$ should
contain a powerful type of ${\rm Th}(\mathcal{A}_i)$ and the
restriction $r(\bar{y})$ of $q(\bar{x})$ to the coordinates
realized by $p_\infty(x)$ should  be powerful for the theory ${\rm
Th}(\mathcal{A}_\infty)$, where $\mathcal{A}_\infty$ is infinite
and saturated, as well as realizing $r(\bar{y})$ in a model
$\mathcal{M}\models T$, all types with coordinates satisfying
$p_\infty(x)$ should be realized in $\mathcal{M}$ too. As shown in
\cite{SuCCMCT, Su071}, the type $r(\bar{y})$ has the local
realizability property and satisfies the following conditions: for
each formula $\varphi(\bar{y})\in r(\bar{y})$, there exists a
formula $\psi(\bar{y},\bar{z})$ of $T$ {\rm (}where
$l(\bar{y})=l(\bar{z})${\rm )}, satisfying the following
conditions:

{\rm (i)} for each $\bar{a}\in r(M)$, the formula
$\psi(\bar{a},\bar{y})$ is equivalent to a disjunction  of
principal  formulas  $\psi_i(\bar{a},\bar{y})$,  $i\leq m$, such
 that  $\psi_i(\bar{a},\bar{y})\vdash r(\bar{y})$, and
$\models\psi_i(\bar{a},\bar{b})$ implies, that $\bar{b}$ does not
semi-isolate~$\bar{a}$;

{\rm (ii)} for  every \ $\bar{a},\bar{b}\in r(M)$,  there  exists
 a  tuple  $\bar{c}$  such  that
$\models\varphi(\bar{c})\wedge\psi(\bar{c},\bar{a})\wedge\psi(\bar{c},\bar{b})$.

Since the type $p_\infty(x)$ is not isolated each formula
$\varphi(\bar{y})\in r(\bar{y})$ has realizations $\bar{d}$ in
$\bigcup\limits_{i\in I}A_i$. On the other hand, as we consider
the disjoint union of $\mathcal{A}_i$ and there are no non-trivial
links between distinct $P_i$ and $P_{i'}$, the sets of solutions
for $\psi(\bar{d},\bar{y})$ with $\models\varphi(\bar{d})$ in
$\{\neg P_i(x)\mid \models P_i(d_j)\mbox{ for some
}d_j\in\bar{d}\}$ are either equal or empty being composed by
definable sets without parameters. If these sets are nonempty the
item (i) can not be satisfied: $\psi(\bar{a},\bar{y})$ is not
equivalent to a disjunction  of  principal  formulas. Otherwise
all $\psi$-links for realizations of $r(\bar{y})$ are situated
inside the set of solutions for
$\bar{p}_\infty(\bar{y})=\bigcup\limits_{y_j\in\bar{y}}p_\infty(y_j)$.
In this case for $\bar{a}\models r(\bar{y})$ the formula
$\exists\bar{z}(\psi(\bar{z},\bar{a})\wedge\psi(\bar{z},\bar{y}))$
does not cover the set $r(M)$ since it does not cover each
$\varphi$-approximation of $r(M)$. Thus, the property (ii) fails.

Hence, (i) and (ii) can not be satisfied, there are no powerful
types, and the theory $T$ is not Ehrenfeucht.~$\Box$

\section{Variations of structures related to combinations and $E$-representability}

Clearly, for a disjoint $P$-combination $\mathcal{A}_P$ with
infinite $I$, there is a structure
$\mathcal{A}'\equiv\mathcal{A}_P$ with a structure
$\mathcal{A}'_\infty$. Since the type $p_\infty(x)$ is nonisolated
(omitted in $\mathcal{A}_P$), the cardinalities for
$\mathcal{A}'_\infty$ are unbounded. Infinite structures
$\mathcal{A}'_\infty$ are not necessary elementary equivalent and
can be both elementary equivalent to some $\mathcal{A}_i$ or not.
For instance, if infinitely many structures $\mathcal{A}_i$
contain unary predicates $Q_0$, say singletons, without unary
predicates $Q_1$ and infinitely many $\mathcal{A}_{i'}$ for $i'\ne
i$ contain $Q_1$, say again singletons, without $Q_0$ then
$\mathcal{A}'_\infty$ can contain $Q_0$ without $Q_1$, $Q_1$
without $Q_0$, or both $Q_0$ and $Q_1$. For the latter case,
$\mathcal{A}'_\infty$ is not elementary equivalent neither
$\mathcal{A}_i$, nor $\mathcal{A}_{i'}$.

A natural question arises:

\medskip
{\bf Question 1.} {\em What can be the number of pairwise
elementary non-equivalent structures $\mathcal{A}'_\infty$}?

\medskip
Considering an $E$-combination $\mathcal{A}_E$ with infinite $I$,
and all structures $\mathcal{A}'\equiv\mathcal{A}_E$, there are
two possibilities: each non-empty {\em
$E$-restriction}\index{$E$-restriction} of $\mathcal{A}'_\infty$,
i.~e. a restriction to some $E$-class, is elementary equivalent to
some $\mathcal{A}_i$, $i\in I$, or some $E$-restriction of
$\mathcal{A}'_\infty$ is not elementary equivalent to all
structures $\mathcal{A}_i$, $i\in I$.

Similarly Question 1 we have:

\medskip
{\bf Question 2.} {\em What can be the number of pairwise
elementary non-equivalent $E$-restrictions of structures
$\mathcal{A}'_\infty$}?

\medskip
{\bf Example 3.1.} Let $\mathcal{A}_P$ be a disjoint
$P$-combination with infinite $I$ and composed by infinite
$\mathcal{A}_i$, $i\in I$, such that $I$ is a disjoint union of
infinite $I_j$, $j\in J$, where $\mathcal{A}_{i_j}$ contains only
unary predicates and unique nonempty unary predicate $Q_j$ being a
singleton. Then $\mathcal{A}'_\infty$ can contain any singleton
$Q_j$ and finitely or infinitely many elements in
$\bigcap\limits_{j\in J}\overline{Q}_j$. Thus, there are
$2^{|J|}\cdot(\lambda+1)$ non-isomorphic $\mathcal{A}'_\infty$,
where $\lambda$ is a least upper bound for cardinalities
$\left|\bigcap\limits_{j\in J}\overline{Q}_j\right|$.

\medskip
For $T={\rm Th}(\mathcal{A}_P)$, we denote by
$I_\infty(T,\lambda)$\index{$I_\infty(T,\lambda)$} the number of
pairwise non-isomorphic structures $\mathcal{A}'_\infty$ having
the cardinality $\lambda$.

Clearly, $I_\infty(T,\lambda)\leq I(T,\lambda)$.

If structures $\mathcal{A}'_\infty$ exist and do not have links
with $\mathcal{A}'_P$ (for instance, for a disjoint
$P$-combination) then $I_\infty(T,\lambda)+1\leq I(T,\lambda)$,
since if models of $T$ are isomorphic then their restrictions to
$p_\infty(x)$ are isomorphic too, and $p_\infty(x)$ can be omitted
producing $\mathcal{A}'_\infty=\varnothing$. Here
$I_\infty(T,\lambda)+1=I(T,\lambda)$ if and only if all $I({\rm
Th}(\mathcal{A}_i),\lambda)=1$ and, moreover, for any
$\left(\bigcup\limits_{i\in I} P_i\right)$-restrictions
$\mathcal{B}_P,\mathcal{B}'_P$ of $\mathcal{B},\mathcal{B}'\models
T$ respectively, where $|B|=|B'|=\lambda$, and their
$P_i$-restrictions $\mathcal{B}_i$, $\mathcal{B}'_i$, there
are isomorphisms $f_i\mbox{\rm : }\mathcal{B}_i\,\mbox{\raisebox{0.7\height}{$\sim$}}\!\!\!\!\!%
\mbox{\raisebox{-0.4\height}{$\to$}}\,\mathcal{B}'_i$ preserving
$P_i$ and with an isomorphism $\bigcup\limits_{i\in I}f_i\mbox{\rm : }\mathcal{B}_P\,\mbox{\raisebox{0.7\height}{$\sim$}}\!\!\!\!\!%
\mbox{\raisebox{-0.4\height}{$\to$}}\,\mathcal{B}'_P$.

The following example illustrates the equality
$I_\infty(T,\lambda)+1=I(T,\lambda)$ with some $I({\rm
Th}(\mathcal{A}_i),\lambda)>1$.

\medskip
{\bf Example 3.2.} Let $P_0$ be a unary predicate containing a
copy of the Ehrenfeucht example \cite{Va} with a dense linear
order $\leq$ and an increasing chain of singletons coding
constants $c_k$, $k\in\omega$; $P_n$, $n\geq 1$, be pairwise
disjoint unary predicates disjoint to $P_0$ such that
$P_1=(-\infty,c'_0)$ $P_{n+2}=[c'_n,c'_{n+1})$, $n\in\omega$, and
$\bigcup\limits_{n\geq 1}P_n$ forms a universe of prime model
(over $\varnothing$) for another copy of the Ehrenfeucht example
with a dense linear order $\leq'$ and an increasing chain of
constants $c'_k$, $k\in\omega$. Now we extend the language
$$\Sigma=\langle\leq,\leq',P_n,\{c_n\},\{c'_n\}\rangle_{n\in\omega}$$
by a bijection $f$ between $P_0=\{a\mid a\leq c_0\mbox{ or
}c_0\leq a\}$ and $\{a'\mid a'\leq' c'_0\mbox{ or }c'_0\leq' a'\}$
such that $a\leq b\Leftrightarrow f(a)\leq' f(b)$. The structures
$\mathcal{A}'_\infty$ consist of realizations $p_\infty(x)$ which
are bijective with realizations of the type $\{c_n<x\mid
n\in\omega\}$.

For the theory $T$ of the described structure ${\rm
EComb}_P(\mathcal{A}_i)_{i\in I}$ we have $I(T,\omega)=3$ (as for
the Ehrenfeucht example and the restriction of $T$ to $P_0$) and
$I_\infty(T,\omega)=2$ (witnessed by countable structures with
least realizations of $p_\infty(x)$ and by countable structure
with realizations of $p_\infty(x)$ all of which are not least).

\medskip
For Example 3.1 of a theory $T$ with singletons $Q_j$ in
$\mathcal{A}_i$ and for a cardinality $\lambda\geq 1$, we have
$$I_\infty(T,\lambda)=\left\{\begin{array}{rl} \sum\limits_{i=0}^{{\rm
min}\{|J|,\lambda\}}C^i_{|J|}, & \mbox{ if } J\mbox{ and }\lambda\mbox{ are finite;} \\
|J|, & \mbox{ if } J \mbox{ is infinite and } |J|>\lambda;
\\
2^{|J|}, & \mbox{ if } J\mbox{ is infinite and } |J|\leq\lambda.
\end{array}\right.$$


\medskip
Clearly, $\mathcal{A}'\equiv\mathcal{A}_P$ realizing $p_\infty(x)$
is not elementary embeddable into $\mathcal{A}_P$ and can not be
represented as a disjoint $P$-combination of
$\mathcal{A}'_i\equiv\mathcal{A}_i$, $i\in I$. At the same time,
there are $E$-combinations such that all
$\mathcal{A}'\equiv\mathcal{A}_E$ can be represented as
$E$-combinations of some $\mathcal{A}'_j\equiv\mathcal{A}_i$. We
call this representability of $\mathcal{A}'$ to be the {\em
$E$-representability}. If, for instance, all $\mathcal{A}_i$ are
infinite structures of the empty language then any
$\mathcal{A}'\equiv\mathcal{A}_E$ is an $E$-combination of some
infinite structures $\mathcal{A}'_j$ of the empty language too.

Thus we have:

\medskip
{\bf Question 3.} {\em What is a characterization of
$E$-representability for all $\mathcal{A}'\equiv\mathcal{A}_E$}?

\medskip
{\bf Definition} (cf. \cite{Henk}). For a first-order formula
$\varphi(x_1,\ldots,x_n)$, an equivalence relation $E$ and a
formula $\sigma(x)$ we define a {\em
$(E,\sigma)$-relativized}\index{Formula!$(E,\sigma)$-relativized}
formula $\varphi^{E,\sigma}$\index{$\varphi^{E,\sigma}$} by
induction:

$({\rm i})$ if $\varphi$ is an atomic formula then
$\varphi^{E,\sigma}=\varphi(x_1,\ldots,x_n)\wedge\bigwedge\limits_{i,j=1}^n
E(x_i,x_j)\wedge\exists y(E(x_1,y)\wedge\sigma(y))$;

$({\rm ii})$ if $\varphi=\psi\tau\chi$, where
$\tau\in\{\wedge,\vee,\to\}$, and $\psi^{E,\sigma}$ and
$\chi^{E,\sigma}$ are defined then
$\varphi^{E,\sigma}=\psi^{E,\sigma}\tau\chi^{E,\sigma}$;

$({\rm iii})$ if
$\varphi(x_1,\ldots,x_n)=\neg\psi(x_1,\ldots,x_n)$ and
$\psi^{E,\sigma}(x_1,\ldots,x_n)$ is defined then
$\varphi^{E,\sigma}(x_1,\ldots,x_n)=\neg\psi^{E,\sigma}(x_1,\ldots,x_n)\wedge\bigwedge\limits_{i,j=1}^n
(E(x_i,x_j)\wedge\exists y(E(x_1,y)\wedge\sigma(y))$;

$({\rm iv})$ if $\varphi(x_1,\ldots,x_n)=\exists
x\psi(x,x_1,\ldots,x_n)$ and $\psi^{E,\sigma}(x,x_1,\ldots,x_n)$
is defined then
$$\varphi^{E,\sigma}(x_1,\ldots,x_n)=\exists
x\left(\bigwedge\limits_{i=1}^n (E(x,x_i)\wedge\exists
y(E(x,y)\wedge\sigma(y))\wedge\psi^{E,\sigma}(x,x_1,\ldots,x_n)\right);$$

$({\rm v})$ if $\varphi(x_1,\ldots,x_n)=\forall
x\psi(x,x_1,\ldots,x_n)$ and $\psi^{E,\sigma}(x,x_1,\ldots,x_n)$
is defined then
$$\varphi^{E,\sigma}(x_1,\ldots,x_n)=\forall
x\left(\bigwedge\limits_{i=1}^n E(x,x_i)\wedge\exists
y(E(x,y)\wedge\sigma(y))\to\psi^{E,\sigma}(x,x_1,\ldots,x_n)\right).$$

\medskip
We write $E$ instead of $(E,\sigma)$ if $\sigma=(x\approx x)$.

Note that two $E$-classes $E_i$ and $E_j$ with structures
$\mathcal{A}_i$ and $\mathcal{A}_j$ (of a language $\Sigma$),
respectively, are not elementary equivalent if and only if there
is a $\Sigma$-sentence $\varphi$ such that
$\mathcal{A}_E\upharpoonright E_i\models\varphi^E$ (with
$\mathcal{A}_i\models\varphi$)  and $\mathcal{A}_E\upharpoonright
E_j\models(\neg\varphi)^E$ (with
$\mathcal{A}_j\models\neg\varphi$). In this case, the formula
$\varphi$ is called {\em
$(i,j)$-separating}\index{Formula!$(i,j)$-separating}.

\medskip
The following properties are obvious:

(1) If $\varphi$ is $(i,j)$-separating then $\neg\varphi$ is
$(j,i)$-separating.

(2) If $\varphi$ is $(i,j)$-separating and $\psi$ is
$(i,k)$-separating then $\varphi\wedge\psi$ is both
$(i,j)$-separating and $(i,k)$-separating.

(3) There is a set $\Phi_i$ of $(i,j)$-separating sentences, for
$j$ in some $J\subseteq I\setminus\{i\}$, which separates
$\mathcal{A}_i$ from all structures
$\mathcal{A}_j\not\equiv\mathcal{A}_i$.

The set $\Phi_i$ is called {\em
$e$-separating}\index{Set!$e$-separating} (for $\mathcal{A}_i$)
and $\mathcal{A}_i$ is {\em
$e$-separable}\index{Structure!$e$-separable} (witnessed by
$\Phi_i$).

Assuming that some $\mathcal{A}'\equiv\mathcal{A}_E$ is not
$E$-representable, we get an $E'$-class with a structure
$\mathcal{B}$ in $\mathcal{A}'$ which is $e$-separable from all
$\mathcal{A}_i$, $i\in I$, by a set $\Phi$. It means that for some
sentences $\varphi_i$ with $\mathcal{A}_E\upharpoonright
E_i\models\varphi_i^E$, i.~e., $\mathcal{A}_i\models\varphi_i$,
the sentences $\left(\bigwedge\limits_{i\in
I_0}\neg\varphi_i\right)^E$, where $I_0\subseteq_{\rm fin}I$, form
a consistent set, satisfying the restriction of $\mathcal{A}'$ to
the class $E'_B$ with the universe $B$ of $\mathcal{B}$.

Thus, answering Question 3 we have

\medskip
{\bf Proposition 3.3.} {\em For any $E$-combination
$\mathcal{A}_E$ the following conditions are equivalent:

$(1)$ there is $\mathcal{A}'\equiv\mathcal{A}_E$ which is not
$E$-representable;

$(2)$ there are sentences $\varphi_i$ such that
$\mathcal{A}_i\models\varphi_i$, $i\in I$, and the set of
sentences $\left(\bigwedge\limits_{i\in
I_0}\neg\varphi_i\right)^E$, where $I_0\subseteq_{\rm fin}I$, is
consistent with ${\rm Th}(\mathcal{A}_E)$.}

\medskip
Proposition 3.3 implies

\medskip
{\bf Corollary 3.4.} {\em If $\mathcal{A}_E$ has only finitely
many pairwise elementary non-equivalent $E$-classes then each
$\mathcal{A}'\equiv\mathcal{A}_E$ is $E$-representable.}

\section{$e$-spectra}

If there is $\mathcal{A}'\equiv\mathcal{A}_E$ which is not
$E$-representable, we have the $E'$-representability replacing $E$
by $E'$ such that $E'$ is obtained from $E$ adding equivalence
classes with models for all theories $T$, where $T$ is a theory of
a restriction $\mathcal{B}$ of a structure
$\mathcal{A}'\equiv\mathcal{A}_E$ to some $E$-class and
$\mathcal{B}$ is not elementary equivalent to the structures
$\mathcal{A}_i$. The resulting structure $\mathcal{A}_{E'}$ (with
the $E'$-representability) is a {\em
$e$-completion}\index{$e$-completion}, or a {\em
$e$-saturation}\index{$e$-saturation}, of $\mathcal{A}_{E}$. The
structure $\mathcal{A}_{E'}$ itself is called {\em
$e$-complete}\index{Structure!$e$-complete}, or {\em
$e$-saturated}\index{Structure!$e$-saturated}, or {\em
$e$-universal}\index{Structure!$e$-universal}, or {\em
$e$-largest}\index{Structure!$e$-largest}.

For a structure $\mathcal{A}_E$ the number of {\em
new}\index{Structure!new} structures with respect to the
structures $\mathcal{A}_i$, i.~e., of the structures $\mathcal{B}$
which are pairwise elementary non-equivalent and elementary
non-equivalent to the structures $\mathcal{A}_i$, is called the
{\em $e$-spectrum}\index{$e$-spectrum} of $\mathcal{A}_E$ and
denoted by $e$-${\rm Sp}(\mathcal{A}_E)$.\index{$e$-${\rm
Sp}(\mathcal{A}_E)$} The value ${\rm sup}\{e$-${\rm
Sp}(\mathcal{A}'))\mid\mathcal{A}'\equiv\mathcal{A}_E\}$ is called
the {\em $e$-spectrum}\index{$e$-spectrum} of the theory ${\rm
Th}(\mathcal{A}_E)$ and denoted by $e$-${\rm Sp}({\rm
Th}(\mathcal{A}_E))$.\index{$e$-${\rm Sp}({\rm
Th}(\mathcal{A}_E))$}

If $\mathcal{A}_E$ does not have $E$-classes $\mathcal{A}_i$,
which can be removed, with all $E$-classes
$\mathcal{A}_j\equiv\mathcal{A}_i$, preserving the theory ${\rm
Th}(\mathcal{A}_E)$, then $\mathcal{A}_E$ is called {\em
$e$-prime}\index{Structure!$e$-prime}, or {\em
$e$-minimal}\index{Structure!$e$-minimal}.

For a structure $\mathcal{A}'\equiv\mathcal{A}_E$ we denote by
${\rm TH}(\mathcal{A}')$ the set of all theories ${\rm
Th}(\mathcal{A}_i)$\index{${\rm Th}(\mathcal{A}_i)$} of
$E$-classes $\mathcal{A}_i$ in $\mathcal{A}'$.

By the definition, an $e$-minimal structure $\mathcal{A}'$
consists of $E$-classes with a minimal set ${\rm
TH}(\mathcal{A}')$. If ${\rm TH}(\mathcal{A}')$ is the least for
models of ${\rm Th}(\mathcal{A}')$ then $\mathcal{A}'$ is called
{\em $e$-least}.\index{Structure!$e$-least}

\medskip
The following proposition is obvious:

\medskip
{\bf Proposition 4.1.} {\em $1.$ For a given language $\Sigma$,
$0\leq e$-${\rm Sp}({\rm Th}(\mathcal{A}_E))\leq 2^{{\rm
max}\{|\Sigma|,\omega\}}$.

\medskip
$2.$ A structure $\mathcal{A}_{E}$ is $e$-largest if and only if
$e$-${\rm Sp}(\mathcal{A}_E)=0$. In particular, an $e$-minimal
structure $\mathcal{A}_E$ is $e$-largest is and only if $e$-${\rm
Sp}({\rm Th}(\mathcal{A}_E))=0$.

\medskip
$3.$ Any weakly saturated structure $\mathcal{A}_E$, i.~e., a
structure realizing all types of ${\rm Th}(\mathcal{A}_E)$ is
$e$-largest.

\medskip
$4.$ For any $E$-combination $\mathcal{A}_E$, if $\lambda\leq
e$-${\rm Sp}({\rm Th}(\mathcal{A}_E))$ then there is a structure
$\mathcal{A}'\equiv\mathcal{A}_E$ with $e$-${\rm
Sp}(\mathcal{A}')=\lambda$; in particular, any theory ${\rm
Th}(\mathcal{A}_E)$ has an $e$-largest model.

\medskip
$5.$ For any structure $\mathcal{A}_{E}$, $e\mbox{-}{\rm
Sp}(\mathcal{A}_E)=|{\rm TH}(\mathcal{A}'_{E'})\setminus{\rm
TH}(\mathcal{A}_{E})|$, where $\mathcal{A}'_{E'}$ is an
$e$-largest model of ${\rm Th}(\mathcal{A}_E)$.

\medskip
$6.$ Any prime structure $\mathcal{A}_E$ is $e$-minimal {\rm (}but
not vice versa as the $e$-minimality is preserved, for instance,
extending an infinite $E$-class of given structure to a greater
cardinality{\rm )}. Any small theory ${\rm Th}(\mathcal{A}_E)$ has
an $e$-minimal model {\rm (}being prime{\rm )}, and in this case,
the structure $\mathcal{A}_E$ is $e$-minimal if and only if
$${\rm
TH}(\mathcal{A}_E)=\bigcap\limits_{\mathcal{A}'\equiv\mathcal{A}_E}{\rm
TH}(\mathcal{A}'),$$ i.~e., $\mathcal{A}_E$ is $e$-least.

\medskip
$7.$ If $\mathcal{A}_E$ is $e$-least then $e$-${\rm
Sp}(\mathcal{A}_E)=e$-${\rm Sp}({\rm Th}(\mathcal{A}_E))$.

\medskip
$8.$ If $e$-${\rm Sp}({\rm Th}(\mathcal{A}_E))$ finite and ${\rm
Th}(\mathcal{A}_E)$ has $e$-least model then $\mathcal{A}_E$ is
$e$-minimal if and only if $\mathcal{A}_E$ is $e$-least and if and
only if $e$-${\rm Sp}(\mathcal{A}_E)=e$-${\rm Sp}({\rm
Th}(\mathcal{A}_E))$.

\medskip
$9.$ If $e$-${\rm Sp}({\rm Th}(\mathcal{A}_E))$ is infinite then
there are $\mathcal{A}'\equiv \mathcal{A}_E$ such that $e$-${\rm
Sp}(\mathcal{A}')=e$-${\rm Sp}({\rm Th}(\mathcal{A}_E))$ but
$\mathcal{A}'$ is not $e$-minimal.

\medskip
$10.$ A countable $e$-minimal structure $\mathcal{A}_E$ is prime
if and only if each $E$-class $\mathcal{A}_i$ is a prime
structure.}

\medskip
Reformulating Proposition 2.2 we have

\medskip
{\bf Proposition 4.2.} {\em For $E$-combinations which are not
{\rm EComb}, a countable theory ${\rm Th}(\mathcal{A}_E)$ without
finite models is $\omega$-categorical if and only if $e$-${\rm
Sp}({\rm Th}(\mathcal{A}_E))=0$ and each $E$-class $\mathcal{A}_i$
is either finite or $\omega$-categorical.}

\medskip
Note that if there are no links between $E$-classes (i.~e., the
Comb is considered, not EComb) and there is
$\mathcal{A}'\equiv\mathcal{A}_E$ which is not $E$-representable,
then by Compactness the $e$-completion can vary adding arbitrary
(finitely or infinitely) many new $E$-classes with a fixed
structure which is not elementary equivalent to structures in old
$E$-classes.

\medskip
{\bf Proposition 4.3.} {\em For any cardinality $\lambda$ there is
a theory $T={\rm Th}(\mathcal{A}_E)$ of a language $\Sigma$ such
that $|\Sigma|=|\lambda+1|$ and $e$-${\rm Sp}(T)=\lambda$.}

\medskip
{\bf\em Proof.} Clearly, for structures $\mathcal{A}_i$ of fixed
cardinality and with empty language we have $e$-${\rm Sp}({\rm
Th}(\mathcal{A}_E))=0$. For $\lambda>0$ we take a language
$\Sigma$ consisting of unary predicate symbols $P_i$, $i<\lambda$.
Let $\mathcal{A}_{i,n+1}$ be a structure having a universe
$A_{i,n}$ with $n$ elements and $P_i=A_{i,n}$, $P_j=\varnothing$,
$i,j<\lambda$, $i\ne j$, $n\in\omega\setminus\{0\}$. Clearly, the
structure $\mathcal{A}_E$, formed by all $\mathcal{A}_{i,n}$, is
$e$-minimal. It produces structures
$\mathcal{A}'\equiv\mathcal{A}_E$ containing $E$-classes with
infinite predicates $P_i$, and structures of these classes are not
elementary equivalent to the structures $\mathcal{A}_{i,n}$. Thus,
for the theory $T={\rm Th}(\mathcal{A}_E)$ we have $e$-${\rm
Sp}(T)=\lambda$.~$\Box$

\medskip
In Proposition 4.3, we have $e$-${\rm Sp}(T)=|\Sigma(T)|$. At the
same time the following proposition holds.

\medskip
{\bf Proposition 4.4.} {\em For any infinite cardinality $\lambda$
there is a theory $T={\rm Th}(\mathcal{A}_E)$ of a language
$\Sigma$ such that $|\Sigma|=\lambda$ and $e$-${\rm
Sp}(T)=2^\lambda$.}

\medskip
{\bf\em Proof.} Let $P_j$ be unary predicate symbols, $j<\lambda$,
forming the language $\Sigma$, and $\mathcal{A}_i$ be structures
consisting of only finitely many nonempty predicates
$P_{j_1},\ldots,P_{j_k}$ and such that these predicates are
independent. Taking for the structures $\mathcal{A}_i$ all
possibilities for cardinalities of sets of solutions for formulas
$P^{\delta_{j_1}}_{j_1}(x)\wedge\ldots\wedge
P^{\delta_{j_k}}_{j_k}(x)$, $\delta_{j_l}\in\{0,1\}$, we get an
$e$-minimal structure $\mathcal{A}_E$ such that for the theory
$T={\rm Th}(\mathcal{A}_E)$ we have $e$-${\rm Sp}(T)=2^\lambda$.

Another approach for $e$-${\rm Sp}(T)=2^\lambda$ was suggested by
E.A.~Palyutin. Taking infinitely many $\mathcal{A}_i$ with
arbitrarily finitely many disjoint singletons
$R_{j_1},\ldots,R_{j_k}$, where $\Sigma$ consists of $R_j$,
$j<\lambda$, we get $\mathcal{A}'\equiv\mathcal{A}_E$ with
arbitrarily many singletons for any subset of $\lambda$ producing
$2^\lambda$ $E$-classes which are pairwise elementary
non-equivalent.~$\Box$

\medskip
If $e$-${\rm Sp}(T)=0$ the theory $T$ is called {\em
$e$-non-abnormalized}\index{Theory!$e$-non-abnormalized} or {\em
$(e,0)$-abnormalized}\index{Theory!$(e,0)$-abnormalized}.
Otherwise, i.~e., if $e$-${\rm Sp}(T)>0$, $T$ is {\em
$e$-abnormalized}\index{Theory!$e$-abnormalized}. An
$e$-abnormalized theory $T$ with $e$-${\rm Sp}(T)=\lambda$ is
called {\em
$(e,\lambda)$-abnormalized}\index{Theory!$(e,\lambda)$-abnormalized}.
In particular, an $(e,1)$-abnormalized theory is {\em
$e$-categorical}\index{Theory!$e$-categorical}, an
$(e,n)$-abnormalized theory with $n\in\omega\setminus\{0,1\}$ is
{\em $e$-Ehrenfeucht}\index{Theory!$e$-Ehrenfeucht}, an
$(e,\omega)$-abnormalized theory is {\em
$e$-countable}\index{Theory!$e$-countable}, and an
$(e,2^\lambda)$-abnormalized theory is {\em
$(e,\lambda)$-maximal}\index{Theory!$(e,\lambda)$-maximal}.

If $e$-${\rm Sp}(T)=\lambda$ and $T$ has a model $\mathcal{A}_E$
with $e$-${\rm Sp}(\mathcal{A}_E)=\mu$ then $\mathcal{A}_E$ is
called {\em
$(e,\varkappa)$-abnormalized}\index{Model!$(e,\varkappa)$-abnormalized},
where $\varkappa$ is the least cardinality with
$\mu+\varkappa=\lambda$.

By proofs of Propositions 4.3 and 4.4 we have

\medskip
{\bf Corollary 4.5.} {\em For any cardinalities $\mu\leq\lambda$
and the least cardinality $\varkappa$ with $\mu+\varkappa=\lambda$
there is an $(e,\lambda)$-abnormalized theory $T$ with an
$(e,\varkappa)$-abnormalized model $\mathcal{A}_E$.}

\medskip
Let $\mathcal{A}_{E}$ and $\mathcal{B}_{E'}$ be structures and
$\mathcal{C}_{E''}=\mathcal{A}_{E}\coprod\mathcal{B}_{E'}$ be
their disjoint union, where $E''=E\coprod E'$. We denote by ${\rm
ComLim}(\mathcal{A}_{E},\mathcal{B}_{E'})$\index{${\rm
ComLim}(\mathcal{A}_{E},\mathcal{B}_{E'})$} the number of
elementary pairwise non-equivalent structures $\mathcal{D}$ which
are both a restriction of $\mathcal{A}'\equiv \mathcal{A}_{E}$ to
some $E$-class and a restriction of $\mathcal{B}'\equiv
\mathcal{B}_{E'}$ to some $E'$-class as well as $\mathcal{D}$ is
not elementary equivalent to the structures $\mathcal{A}_i$ and
$\mathcal{B}_j$.

\medskip
We have:
$$
{\rm ComLim}(\mathcal{A}_{E},\mathcal{B}_{E'})\leq{\rm
min}\{e\mbox{-}{\rm Sp}({\rm Th}(\mathcal{A}_{E})),e\mbox{-}{\rm
Sp}({\rm Th}(\mathcal{B}_{E'}))\},
$$
$$
{\rm max}\{e\mbox{-}{\rm Sp}({\rm
Th}(\mathcal{A}_{E})),e\mbox{-}{\rm Sp}({\rm
Th}(\mathcal{B}_{E'}))\}\leq e\mbox{-}{\rm Sp}({\rm
Th}(\mathcal{C}_{E''})),
$$
$$
e\mbox{-}{\rm Sp}({\rm Th}(\mathcal{A}_{E}))+e\mbox{-}{\rm
Sp}({\rm Th}(\mathcal{B}_{E'}))=e\mbox{-}{\rm Sp}({\rm
Th}(\mathcal{C}_{E''}))+{\rm
ComLim}(\mathcal{A}_{E},\mathcal{B}_{E'}).
$$

Indeed, all structures witnessing the value $e\mbox{-}{\rm
Sp}({\rm Th}(\mathcal{C}_{E''}))$ can be obtained by ${\rm
Th}(\mathcal{A}_{E})$ or ${\rm Th}(\mathcal{B}_{E'})$ and common
structures are counted for ${\rm
ComLim}(\mathcal{A}_{E},\mathcal{B}_{E'})$.

\medskip
If $\mathcal{A}_{E}=\mathcal{B}_{E'}$ then ${\rm
ComLim}(\mathcal{A}_{E},\mathcal{B}_{E'})=e$-${\rm Sp}({\rm
Th}(\mathcal{A}_{E}))$. Assuming that $\mathcal{A}_{E}$ and
$\mathcal{B}_{E'}$ do not have elementary equivalent classes
$\mathcal{A}_i$ and $\mathcal{B}_j$, the number ${\rm
ComLim}(\mathcal{A}_{E},\mathcal{B}_{E'})$ can vary from $0$ to
$2^{|\Sigma|+\omega}$.

Indeed, if ${\rm Th}(\mathcal{A}_{E})$ or ${\rm
Th}(\mathcal{B}_{E'})$ does not produce new, elementary
non-equivalent classes then ${\rm
ComLim}(\mathcal{A}_{E},\mathcal{B}_{E'})=0$. Otherwise we can
take structures $\mathcal{A}_i$ and $\mathcal{B}_i$ with one unary
predicate symbol $P$ such that $P$ has $2i$ elements for
$\mathcal{A}_i$ and $2i+1$ elements for $\mathcal{B}_i$,
$i\in\omega$. In this case we have ${\rm Sp}({\rm
Th}(\mathcal{A}_{E}))=1$, ${\rm Sp}({\rm
Th}(\mathcal{B}_{E'}))=1$, ${\rm
ComLim}(\mathcal{A}_{E},\mathcal{B}_{E'})=1$, and
$\mathcal{C}_{E''}$ witnessed by structures with infinite
interpretations for $P$. Extending the language by unary
predicates $P_i$, $i<\lambda$, and interpreting $P_i$ in disjoint
structures as for $P$ above, we get ${\rm Sp}({\rm
Th}(\mathcal{A}_{E}))=\lambda$, ${\rm Sp}({\rm
Th}(\mathcal{B}_{E'}))=\lambda$, ${\rm
ComLim}(\mathcal{A}_{E},\mathcal{B}_{E'})=\lambda$. Thus we have

\medskip
{\bf Proposition 4.6.} {\em For any cardinality $\lambda$ there
are structures $\mathcal{A}_E$ and $\mathcal{B}_{E'}$ of a
language $\Sigma$ such that $|\Sigma|=|\lambda+1|$ and ${\rm
ComLim}(\mathcal{A}_{E},\mathcal{B}_{E'})=\lambda$.}

\medskip
Applying proof of Proposition 4.4 with even and odd cardinalities
for intersections of predicates in $\mathcal{A}_i$ and
$\mathcal{B}_j$ respectively, we have ${\rm Sp}({\rm
Th}(\mathcal{A}_{E}))=2^\lambda$, ${\rm Sp}({\rm
Th}(\mathcal{B}_{E'}))=2^\lambda$, ${\rm
ComLim}(\mathcal{A}_{E},\mathcal{B}_{E'})=2^\lambda$. In
particular, we get

\medskip
{\bf Proposition 4.7.} {\em For any infinite cardinality $\lambda$
are structures $\mathcal{A}_E$ and $\mathcal{B}_{E'}$ of a
language $\Sigma$ such that $|\Sigma|=\lambda$ and ${\rm
ComLim}(\mathcal{A}_{E},\mathcal{B}_{E'})=2^\lambda$.}

\medskip
Replacing $E$-classes by unary predicates $P_i$ (not necessary
disjoint) being universes for structures $\mathcal{A}_i$ and
restricting models of ${\rm Th}(\mathcal{A}_P)$ to the set of
realizations of $p_\infty(x)$ we get the {\em
$e$-spectrum}\index{$e$-spectrum} $e$-${\rm Sp}({\rm
Th}(\mathcal{A}_P))$\index{$e$-${\rm Sp}({\rm
Th}(\mathcal{A}_P))$}, i.~e., the number of pairwise elementary
non-equivalent restrictions of $\mathcal{M}\models{\rm
Th}(\mathcal{A}_P)$ to $p_\infty(x)$. We also get the notions of
$(e,\lambda)$-abnormalized theory ${\rm Th}(\mathcal{A}_P)$, of
$(e,\lambda)$-abnormalized model of ${\rm Th}(\mathcal{A}_P)$, and
related notions.

Note that for any countable theory $T={\rm Th}(\mathcal{A}_P)$,
$e$-${\rm Sp}(T)\leq I(T,\omega)$. In particular, if $I(T,\omega)$
is finite then $e$-${\rm Sp}(T)$ is finite too. Moreover, if $T$
is $\omega$-categorical then $e$-${\rm Sp}(T)=0$, and if $T$ is an
Ehrenfeucht theory, then $e$-${\rm Sp}(T)<I(T,\omega)$.
Illustrating the finiteness for Ehrenfeucht theories we consider

\medskip
{\bf Example 4.8.} Similar to Example 3.2, let $T_0$ be the
Ehrenfeucht theory of a structure ${\cal M}_0$, formed from the
structure $\langle\mathbb Q;<\rangle$ by adding singletons $R_k$
for elements $c_k$, $c_k<c_{k+1}$, $k\in\omega$, such that
$\lim\limits_{k\to\infty}c_k=\infty$. It is well known that the
theory $T_3$ has exactly $3$ pairwise non-isomorphic models:

(a) a prime model ${\cal M}_0$
($\lim\limits_{k\to\infty}c_k=\infty$);\index{${\cal M}^n$}

(b) a prime model ${\cal M}_1$ over a realization of powerful type
$p_\infty(x)\in S^1(\varnothing)$, isolated \ by \ sets \ of \
formulas \ $\{c_k<x\mid k\in\omega\}$;

(c) a saturated model ${\cal M}_2$ (the limit
$\lim\limits_{k\to\infty}c_k$ is irrational).

Now we introduce unary predicates $P_i=\{a\in M_0\mid a<c_i\}$,
$i<\omega$, on ${\cal M}_0$. The structures
$\mathcal{A}_i=\mathcal{M}_0\upharpoonright P_i$ form the
$P$-combination $\mathcal{A}_P$ with the universe $M_0$.
Realizations of the type $p_\infty(x)$ in ${\cal M}_1$ and in
$\mathcal{M}_2$ form two elementary non-equivalent structures
$\mathcal{A}_\infty$ and $\mathcal{A}'_\infty$ respectively, where
$\mathcal{A}_\infty$ has a dense linear order with a least element
and $\mathcal{A}'_\infty$ has a dense linear order without
endpoints. Thus, $e$-${\rm Sp}(T_0)=2$ and $T_0$ is
$e$-Ehrenfeucht.

As E.A.~Palyutin noticed, varying unary predicates $P_i$ in the
following way: $P_{2i}=\{a\in M_0\mid a<c_{2i}\}$,
$P_{2i+1}=\{a\in M_0\mid a\leq c_{2i+1}\}$, we get $e$-${\rm
Sp}(T_3)=4$ since the structures $\mathcal{A}'_\infty$ have dense
linear orders with(out) least elements and with(out) greatest
elements.

Modifying Example above, let $T_n$ be the Ehrenfeucht theory of a
structure ${\cal M}^n$, formed from the structure $\langle\mathbb
Q;<\rangle$ by adding constants $c_k$, $c_k<c_{k+1}$,
$k\in\omega$, such that $\lim\limits_{k\to\infty}c_k=\infty$, and
unary predicates $R_0,\ldots,R_{n-2}$ which form a~partition of
the set $\mathbb Q$ of rationals, with
$$\models\forall x,y\:((x<y)\to\exists z\:((x<z)\wedge (z<y)\wedge
R_i(z))),\mbox{\ } i=0,\ldots,n-2.$$ The theory $T_n$ has exactly
$n+1$ pairwise non-isomorphic models:

(a) a prime model ${\cal M}^n$
($\lim\limits_{k\to\infty}c_k=\infty$);

(b) prime models ${\cal M}_i^n$ over realizations of powerful
types $p_i(x)\in S^1(\varnothing)$, isolated \ by \ sets \ of \
formulas \ $\{c_k<x\mid k\in\omega\}\cup\{P_i(x)\}$,
$i=$~$0,\ldots,n-2$ ($\lim\limits_{k\to\infty}c_k\in P_i$);

(c) a saturated model ${\cal M}_infty^n$ (the limit
$\lim\limits_{k\to\infty}c_k$ is irrational).

Now we introduce unary predicates $P_i=\{a\in M^n\mid a<c_i\}$,
$i<\omega$, on ${\cal M}^n$. The structures
$\mathcal{A}_i=\mathcal{M}^n\upharpoonright P_i$ form the
$P$-combination $\mathcal{A}_P$ with the universe $M^n$.
Realizations of the type $p_\infty(x)$ in ${\cal M}_i^n$ and in
$\mathcal{M}_\infty^n$ form $n-1$ elementary non-equivalent
structures $\mathcal{A}^n_j$, $j\leq n-2$, and
$\mathcal{A}^n_{\infty}$, where $\mathcal{A}^n_j$ has a dense
linear order with a least element in $R_j$, and
$\mathcal{A}^n_\infty$ has a dense linear order without endpoints.
Thus, $e$-${\rm Sp}(T_n)=n$ and $T_n$ is $e$-Ehrenfeucht.

Note that in the example above the type $p_\infty(x)$ has $n-1$
completions by formulas $R_0(x),\ldots,R_{n-2}(x)$.

\medskip
{\bf Example 4.9.} Taking a disjoint union $\mathcal{M}$ of
$m\in\omega\setminus\{0\}$ copies of $\mathcal{M}_0$ in the
language $\{<_j, R_k\}_{j<m,k\in\omega}$ and unary predicates
$P_i=\{a\mid\mathcal{M}\models\exists x(a<x\wedge R_i(x))\}$ we
get the $P$-combination $\mathcal{A}_P$ with the universe $M$ for
the structures $\mathcal{A}_i=\mathcal{M}\upharpoonright P_i$,
$i\in\omega$. We have $e$-${\rm Sp}({\rm
Th}(\mathcal{A}_P))=3^m-1$ since each connected component of
$\mathcal{M}$ produces at most two possibilities for dense linear
orders or can be empty on the set of realizations of
$p_\infty(x)$, and at least one connected component has
realizations of $p_\infty(x)$.

Marking the relations $<_j$ by the same symbol $<$ we get the
theory $T$ with $$e\mbox{-}{\rm
Sp}(T)=\sum\limits_{l=1}^m(l+1)=\frac{m(m+1)}{2}+m=\frac{m^2+3m}{2}.$$

\medskip
Examples 4.8 and 4.9 illustrate that having a powerful type
$p_\infty(x)$ we get $e$-${\rm Sp}({\rm Th}(\mathcal{A}_P))\ne 1$,
i.~e., there are no $e$-categorical theories ${\rm
Th}(\mathcal{A}_P)$ with a powerful type $p_\infty(x)$. Moreover,
we have

\medskip
{\bf Theorem 4.10.} {\em For any theory ${\rm Th}(\mathcal{A}_P)$
with non-symmetric or definable semi-isolation on the complete
type $p_\infty(x)$, $e$-${\rm Sp}({\rm Th}(\mathcal{A}_P))\ne 1$.}

\medskip
{\bf\em Proof.}
Assuming the hypothesis we take a realization $a$ of $p_\infty(x)$
and construct step-by-step a {\em $(a,p_\infty(x))$-thrifty model}
$\mathcal{N}$ of ${\rm Th}(\mathcal{A}_P)$, i.~e., a model
satisfying the following condition: if $\varphi(x,y)$ is a formula
such that $\varphi(a,y)$ is consistent and there are no consistent
formulas $\psi(a,y)$ with $\psi(a,y)\vdash p_\infty(x)$ then
$\varphi(a,\mathcal{N})=\varnothing$.

At the same time, since $p_\infty(x)$ is non-isolated, for any
realization $a$ of $p_\infty(x)$ the set
$p_\infty(x)\cup\{\neg\varphi(a,x)\mid\varphi(a,x)\vdash
p_\infty(x)\}$ is consistent. Then there is a model
$\mathcal{N}'\models{\rm Th}(\mathcal{A}_P)$ realizing
$p_\infty(x)$ and which is not $(a',p_\infty(x))$-thrifty for any
realization $a'$ of $p_\infty(x)$.

If semi-isolation is non-symmetric, $\mathcal{N}\upharpoonright
p_\infty(x)$ and $\mathcal{N}'\upharpoonright p_\infty(x)$ are not
elementary equivalent since the formula $\varphi(a,y)$ witnessing
the non-symmetry of semi-isolation has solutions in
$\mathcal{N}'\upharpoonright p_\infty(x)$ and does not have
solutions in $\mathcal{N}\upharpoonright p_\infty(x)$.

If semi-isolation is definable and witnessed by a formula
$\psi(a,y)$ then again $\mathcal{N}\upharpoonright p_\infty(x)$
and $\mathcal{N}'\upharpoonright p_\infty(x)$ are not elementary
equivalent since $\neg\psi(a,y)$ is realized in
$\mathcal{N}'\upharpoonright p_\infty(x)$ and it does not have
solutions in $\mathcal{N}\upharpoonright p_\infty(x)$

Thus, $e$-${\rm Sp}({\rm Th}(\mathcal{A}_P))> 1$.~$\Box$

\medskip
Since non-definable semi-isolation implies that there are
infinitely many $2$-types, we have

\medskip
{\bf Corollary 4.11.} {\em For any theory ${\rm
Th}(\mathcal{A}_P)$ with $e$-${\rm Sp}({\rm Th}(\mathcal{A}_P))=1$
the structures $\mathcal{A}'_\infty$ are not
$\omega$-categorical.}

\medskip
Applying modifications of the Ehrenfeucht example as well as
constructions in \cite{SuCCMCT}, the results for $e$-spectra of
$E$-combinations are modified for $P$-combinations:

\medskip
{\bf Proposition 4.12.} {\em For any cardinality $\lambda$ there
is a theory $T={\rm Th}(\mathcal{A}_P)$ of a language $\Sigma$
such that $|\Sigma|={\rm max}\{\lambda,\omega\}$ and $e$-${\rm
Sp}(T)=\lambda$.}

\medskip
{\bf\em Proof.} Clearly, if $p_\infty(x)$ is inconsistent then
$e$-${\rm Sp}(T)=0$. Thus, the assertion holds for $\lambda=0$.

If $\lambda=1$ we take a theory $T_1$ with disjoint unary
predicates $P_i$, $i\in\omega$, and a symmetric irreflexive binary
relation $R$ such that each vertex has $R$-degree $2$, each $P_i$
has infinitely many connected components, and each connected
component on $P_i$ has diameter $i$. Now structures on
$p_\infty(x)$ have connected components of infinite diameter, all
these structures are elementary equivalent, and $e$-${\rm
Sp}(T_1)=1$.

If $\lambda=n>1$ is finite, we take the theory $T_n$ in Example
4.8 with $e$-${\rm Sp}(T_n)=n$, as well as we can take a generic
Ehrenfeucht theory $T'_\lambda$ with ${\rm RK}(T'_\lambda)=2$ and
with $\lambda-1$ limit model $\mathcal{M}_i$ over the type
$p_\infty(x)$, $i<\lambda-1$, such that each $\mathcal{M}_i$ has a
$Q_j$-chains, $j\leq i$, and does not have $Q_k$-chains for $k>i$.
Restricting the limit models to $p_\infty(x)$ we get $\lambda$
elementary non-equivalent structures including the prime structure
$\mathcal{N}^0$ without $Q_i$-chains and structures
$\mathcal{M}_i\upharpoonright p_\infty(x)$, $i<\lambda-1$, which
are elementary non-equivalent by distinct (non)existence of
$Q_j$-chains.

Similarly, taking $\lambda\geq\omega$ disjoint binary predicates
$R_j$ for the Ehrenfeucht example in 4.8 we have $\lambda$
structures with least elements in $R_j$ which are not elementary
equivalent each other. Producing the theory $T_\lambda$ we have
$e$-${\rm Sp}(T_\lambda)=\lambda$.

Modifying the generic Ehrenfeucht example taking $\lambda$ binary
predicates $Q_j$ with $Q_j$-chains which do not imply $Q_k$-chains
for $k>i$ we get $\lambda$ elementary non-equivalent restrictions
to $p_\infty(x)$.~$\Box$

\medskip
Note that as in Example 4.8 the type $p_\infty(x)$ for the
Ehrenfeucht-like example $T_\lambda$ has $\lambda$ completions by
the formulas $R_j(x)$ whereas the type $p_\infty(x)$ for the
generic Ehrenfeucht theory  is complete. At the same time having
$\lambda$ completions for the $p_\infty(x)$-restrictions related
to $T_\lambda$, the $p_\infty(x)$-restrictions the generic
Ehrenfeucht examples with complete $p_\infty(x)$ can violet the
uniqueness of the complete $1$-type like the Ehrenfeucht example
$T_0$, where $\mathcal{A}_\infty$ realizes two complete 1-types:
the type of the least element and the type of elements which are
not least.

\medskip
{\bf Proposition 4.13.} {\em For any infinite cardinality
$\lambda$ there is a theory $T={\rm Th}(\mathcal{A}_P)$ of a
language $\Sigma$ such that $|\Sigma|=\lambda$ and $e$-${\rm
Sp}(T)=2^\lambda$.}

\medskip
{\bf\em Proof.} Let $T$ be the theory of independent unary
predicates $R_j$, $j<\lambda$, (defined by the set of axioms $
\exists x\,(R_{k_1}(x)\wedge\ldots\wedge R_{k_m}(x)\wedge\neg
R_{l_1}(x)\wedge\ldots\wedge\neg R_{l_n}(x)),$ where
$\{k_1,\ldots,k_m\}\cap\{l_1,\ldots,l_n\}=\varnothing$) such that
countably many of them form predicates $P_i$, $i<\omega$, and
infinitely many of them are independent with $P_i$. Thus, $T$ can
be considered as ${\rm Th}(\mathcal{A}_P)$. Restrictions of models
of $T$ to sets of realizations of the type $p_\infty(x)$ witness
that predicates $R_j$ distinct with all $P_i$ are independent.
Denote indexes of these predicates $R_j$ by $J$. Since
$p_\infty(x)$ is non-isolated, for any family
$\Delta=(\delta_j)_{j\in J}$, where $\delta_j\in\{0,1\}$, the
types $q_\Delta(x)=\{R^{\delta_j}_j\mid j\in J\}$ can be pairwise
independently realized and omitted in structures
$\mathcal{M}\upharpoonright p_\infty(x)$ for $\mathcal{M}\models
T$. Then any predicate $R_j$ can be independently realized and
omitted in these restrictions. Thus there are $2^\lambda$
restrictions with distinct theories, i.~e., $e$-${\rm
Sp}(T)=2^\lambda$.~$\Box$

\medskip
Since for $E$-combinations $\mathcal{A}_E$ and $P$-combinations
$\mathcal{A}_P$ and their limit structures $\mathcal{A}_\infty$,
being respectively structures on $E$-classes and $p_\infty(x)$,
the theories ${\rm Th}(\mathcal{A}_\infty)$ are defined by types
restricted to $E(x,y)$ and $p_\infty(x)$, and for any countable
theory there are either countably many types or continuum many
types, Propositions 4.3, 4.4, 4.12, and 4.13 implies the following

\medskip
{\bf Theorem 4.14.} {\em If $T={\rm Th}(\mathcal{A}_E)$ {\rm
(}respectively, $T={\rm Th}(\mathcal{A}_P)${\rm )} is a countable
theory then $e$-${\rm Sp}(T)\in\omega\cup\{\omega,2^\omega\}$. All
values in $\omega\cup\{\omega,2^\omega\}$ have realizations in the
class of countable theories of $E$-combinations {\rm (}of
$P$-combinations{\rm )}.}

\section{Ehrenfeuchtness for $E$-combinations}

\medskip
{\bf Theorem 5.1.} {\em If the language $\bigcup\limits_{i\in
I}\Sigma(\mathcal{A}_i)$ is at most countable and the structure
$\mathcal{A}_E$ is infinite then the theory $T={\rm
Th}(\mathcal{A}_E)$ is Ehrenfeucht if and only if $e$-${\rm
Sp}(T)<\omega$ {\rm (}which is equivalent here to $e$-${\rm
Sp}(T)=0${\rm )} and for an $e$-largest model
$\mathcal{A}_{E'}\models T$ consisting of $E'$-classes
$\mathcal{A}_j$, $j\in J$, the following conditions hold:

${\rm (a)}$ for any $j\in J$, $I({\rm
Th}(\mathcal{A}_j),\omega)<\omega$;

${\rm (b)}$ there are positively and finitely many $j\in J$ such
that $I({\rm Th}(\mathcal{A}_j),\omega)>1$;


${\rm (c)}$ if 
$I({\rm Th}(\mathcal{A}_j),\omega)\leq 1$ then there are always
finitely many $\mathcal{A}_{j'}\equiv\mathcal{A}_j$ or always
infinitely many $\mathcal{A}_{j'}\equiv\mathcal{A}_j$ independent
of $\mathcal{A}_{E'}\models T$.}

\medskip
{\bf\em Proof.} If $e$-${\rm Sp}(T)<\omega$ and the conditions
(a)--(c) hold then the theory $T$ is Ehrenfeucht since each
countable model $\mathcal{A}_{E''}\models T$ is composed of
disjoint models with universes $E''_k=A_k$, $k\in K$, and
$I(T,\omega)$ is a sum $\sum\limits_{l=0}^{e\mbox{-}{\rm Sp}(T)}$
of finitely many possibilities for models with $l$ representatives
with respect to the elementary equivalence of $E''$-classes that
are not presented in a prime (i.~e., $e$-minimal) model of $T$.
These possibilities are composed by finitely many possibilities of
$I({\rm Th}(\mathcal{A}_k),\omega)>1$ for $\mathcal{A}_{k'}\equiv
\mathcal{A}_k$ and finitely many of $\mathcal{A}_{k''}\not\equiv
\mathcal{A}_k$ with $I({\rm Th}(\mathcal{A}_{k''}),\omega)>1$.
Moreover, there are $\hat{C}(I({\rm
Th}(\mathcal{A}_k),\omega),m_i)$ possibilities for substructures
consisting of $\mathcal{A}_{k'}\equiv\mathcal{A}_k$ 
where $m_i$ is the number of $E$-classes having the theory ${\rm
Th}(\mathcal{A}_k)$, $\hat{C}(n,m)=C^{m}_{n+m-1}$ is the number of
combinations with repetitions for $n$-element sets with $m$
places. The formula for $I(T,\omega)$ is based on the property
that each $E''$-class with the structure $\mathcal{A}_k$ can be
replaced, preserving the elementary equivalence of
$\mathcal{A}_{E''}$, by arbitrary
$\mathcal{B}\equiv\mathcal{A}_k$.

Now we assume that the theory $T$ is Ehrenfeucht. Since models of
$T$ with distinct theories of $E$-classes are not isomorphic, we
have $e$-${\rm Sp}(T)<\omega$. Applying the formula for
$I(T,\omega)$ we have the conditions (a), (b). The condition (c)
holds since varying unboundedly many
$\mathcal{A}_{j'}\equiv\mathcal{A}_j$ we get
$I(T,\omega)\geq\omega$.

The conditions $e$-${\rm Sp}(T)<\omega$ and $e$-${\rm Sp}(T)=0$
are equivalent. Indeed, if $e$-${\rm Sp}(T)>0$ then taking an
$e$-minimal model $\mathcal{M}$ we get, by Compactness,
unboundedly many $E$-classes, which are elementary non-equivalent
to $E$-classes in $\mathcal{M}$. It implies that
$I(T,\omega)\geq\omega$.~$\Box$

\medskip
Since any prime structure is $e$-minimal (but not vice versa as
the $e$-minimality is preserved, for instance, extending an
infinite $E$-class of given structure to a greater cardinality
preserving the elementary equivalence) and any Ehrenfeucht theory
$T$, being small, has a prime model, any Ehrenfeucht theory ${\rm
Th}(\mathcal{A}_E)$ has an $e$-minimal model.

\noindent Sobolev Institute of Mathematics, \\ 4, Acad. Koptyug
avenue, Novosibirsk, 630090, Russia; \\ Novosibirsk State
Technical
University, \\ 20, K.Marx avenue, Novosibirsk, 630073, Russia; \\
Novosibirsk State University, \\ 2, Pirogova street, Novosibirsk,
630090, Russia;
\\ Institute of Mathematics and Mathematical Modeling, \\
125, Pushkina Street, Almaty, 050010, Kazakhstan

\end{document}